\documentclass[a4paper,12pt]{amsart}
\usepackage[a4paper,hscale=0.7,vscale=0.75,centering]{geometry}
\usepackage{fullpage}
\usepackage{url}
\usepackage{xcolor}

\newtheorem{theorem}{Theorem}[section]
\newtheorem{lemma}[theorem]{Lemma}
\newtheorem{corollary}[theorem]{Corollary}

\theoremstyle{definition}
\newtheorem{definition}[theorem]{Definition}
\newtheorem{example}[theorem]{Example}

\newtheorem{assumption}{Assumption}
\newtheorem*{proofDragomir}{Proof of Lemma \ref{lemma:Dragomir}}
\newtheorem*{proofExistence}{Proof of Theorem \ref{thm:existence}}
\newtheorem*{proofDifferentiability}{Proof of Theorem \ref{thm:differentiability}}
\newtheorem*{proofPLI}{Proof of Theorem \ref{thm:PLI}}

\renewcommand{\div}{\textrm{div}}

\theoremstyle{remark}
\newtheorem{remark}[theorem]{Remark}

\numberwithin{equation}{section}

\begin{document}

\title{Polyak-\L ojasiewicz inequality on the space of measures and convergence of mean-field birth-death processes}

\author{Linshan Liu}
\address{School of Mathematical and Computer Sciences, Heriot-Watt University, Edinburgh, EH14 4AS, UK}
\email{ll2018@hw.ac.uk}

\author{Mateusz B. Majka}
\address{School of Mathematical and Computer Sciences, Heriot-Watt University, Edinburgh, EH14 4AS, UK}
\email{m.majka@hw.ac.uk}

\author{\L ukasz Szpruch}
\address{School of Mathematics, University of Edinburgh,  James Clerk Maxwell Building, Peter Guthrie Tait Road,
	Edinburgh EH9 3FD, UK}
\email{l.szpruch@ed.ac.uk}



\dedicatory{}

\begin{abstract}
The Polyak-\L ojasiewicz inequality (P\L I) in $\mathbb{R}^d$ is a natural condition for proving convergence of gradient descent algorithms \cite{Karimi2016}. In the present paper, we study an analogue of P\L I on the space of probability measures $\mathcal{P}(\mathbb{R}^d)$ and show that it is a natural condition for showing exponential convergence of a class of birth-death processes related to certain mean-field optimization problems. We verify P\L I for a broad class of such problems for energy functions regularised by the KL-divergence. 
\end{abstract}

\maketitle

\section{Introduction}\label{section:introduction}

Consider a classical optimization problem, where one is interested in finding a global minimum of a differentiable function $f : \mathbb{R}^d \to \mathbb{R}$. A natural condition on $f$, under which the gradient descent  algorithm has a geometric convergence rate to $\min_{y\in \mathbb R^d} f(y)$, is the Polyak-\L ojasiewicz inequality (P\L I) 
\begin{equation}\label{eq:PLIonRd}
\frac{1}{\kappa} \| \nabla f(x)\|^2 \geq f(x) - \min_{y\in \mathbb R^d} f(y) \,,
\end{equation} 
required to hold with a positive constant $\kappa > 0$, for all $x \in \mathbb{R}^d$ (see \cite{Karimi2016} and the references therein, or \cite{Bolte2010, BlanchetBolte2018} for other variants of \L ojasiewicz inequalities). 
It is easy to see that when $f$ is strictly convex, \eqref{eq:PLIonRd} holds, but the converse is not necessarily true. 

In the present paper we are concerned with an optimization problem on the space of probability measures $\mathcal{P}(\mathbb{R}^d)$. We consider a function $V: \mathcal{P}(\mathbb{R}^d) \to \mathbb{R}$, and we want to find a minimizing measure $m^* \in \mathcal{P}(\mathbb{R}^d)$. Such optimization problems have attracted considerable attention in recent years, see e.g.\ \cite{ChizatBach, Montanari, HuRenSiskaSzpruch2021, RenWang2022, Chizat2021}. In this setting, there exist multiple different choices of flows of probability measures $(m_t)_{t \geq 0}$ that can serve as analogues of the gradient descent algorithm in $\mathbb{R}^d$, as well as multiple different choices of conditions on $V$ analogous to \eqref{eq:PLIonRd} that can be used to prove convergence of such flows.

The main example of $V$ considered in this paper is an energy function regularised by the KL-divergence. Consider $F: \mathcal{P}(\mathbb{R}^d) \to \mathbb{R}$ (which can be non-linear) and a probability measure $\pi(dx) \propto e^{-U(x)}dx$ with a potential $U: \mathbb{R}^d \to \mathbb{R}$. For any $\sigma \geq 0$, we put
\begin{equation}\label{eq:Vsigma}
V^{\sigma}(m) = F(m) + \frac{\sigma^2}{2} \operatorname{KL}(m|\pi) \,, \quad m\in \mathcal{P}(\mathbb{R}^d)\,,
\end{equation}
where for any $m \in \mathcal{P}(\mathbb{R}^d)$, 
\begin{equation*}
\operatorname{KL}(m|\pi) =\begin{cases}
	 \int_{\mathbb{R}^d} \log{\left( \frac{m(x)}{\pi(x)}\right)} m(x) dx \, &\text{$m$ absolutely continuous with respect to $\pi$,}  \\
	\infty  &\text{otherwise.}
\end{cases}
\end{equation*} 
It is known (see e.g.\ Proposition 2.5 in \cite{HuRenSiskaSzpruch2021}) that $V^{\sigma}$ is minimized by a measure $m^{\sigma,*} \in \mathcal{P}(\mathbb{R}^d)$ satisfying
\begin{equation}\label{defmstar}
m^{\sigma,*}(x) = \frac{1}{Z} \exp{\left( -\frac{2}{\sigma^2} \left( \frac{\delta F}{\delta m} (m^{\sigma,*},x) + U(x) \right)\right)} \,,
\end{equation}
where $Z$ is the normalising constant, and for any $m \in \mathcal{P}(\mathbb{R}^d)$ and $x \in \mathbb{R}^d$, by $\frac{\delta F}{\delta m} (m,x)$ we denote the flat derivative of $F$ with respect to $m$, in the direction of $x \in \mathbb{R}^d$, evaluated at $m$. For any $m$, $m' \in \mathcal{P}(\mathbb{R}^d)$, the function  $\frac{\delta F}{\delta m} : \mathcal{P}(\mathbb{R}^d) \times \mathbb{R}^d \to \mathbb{R}$ satisfies
\begin{equation*}
F(m') - F(m) = \int_0^1 \int_{\mathbb{R}^d} \frac{\delta F}{\delta m} \left( m + \lambda (m'-m) , x \right) (m'-m)(dx) d\lambda \,.
\end{equation*}
See Appendix~\ref{section:appendix2} for more details on flat derivatives. This notion of derivative appears in the literature under several different names, including the linear functional derivative (see Section 5.4.1 in \cite{CarmonaDelarue2018}) or the first variation \cite{AmbrosioGigliSavare2008}. It is important to note that $\frac{\delta F}{\delta m}$ is defined only up to a constant, i.e., for any $C \in \mathbb{R}$, the function $\frac{\delta F}{\delta m} + C$ is also a flat derivative of $F$. Everywhere in this paper we will adopt a normalizing convention requiring $\int_{\mathbb{R}^d} \frac{\delta F}{\delta m} (m,x) m(dx) = 0$, which then makes the choice of the constant unique.

The objective of this work is to identify a flow of measures $(m_t)_{t\geq 0}$ such that $V^{\sigma}(m_t) \to V^{\sigma}(m^{\sigma,*})$ as $t \to \infty$, as well as conditions that ensure that this convergence is exponential. To this end, we equip the space $\mathcal{P}(\mathbb{R}^d)$ with a suitable distance function $d:\mathcal{P}(\mathbb{R}^d)\times \mathcal{P}(\mathbb{R}^d) \rightarrow \mathbb R$ and consider a corresponding gradient flow, where the form of the flow is dictated by the choice of $d$. 
Our main focus is on the Fisher-Rao metric.

\subsection*{Fisher-Rao Gradient Flow}

Let $\mathcal{P}_{ac}(\mathbb{R}^d)$ be the space of probability measures on $\mathbb{R}^d$ that are absolutely continuous with respect to the Lebesgue measure. Then the Fisher-Rao distance between $\mu_0$, $\mu_1 \in \mathcal{P}_{ac}(\mathbb{R}^d)$ is defined by
\begin{equation*}
\operatorname{FR}(\mu_0,\mu_1) =   \int_{\mathbb{R}^d} \left|\sqrt{\mu_0(x)} - \sqrt{\mu_1(x)}  \right|^2 dx \,.
\end{equation*}
One can also consider a dynamic representation of the Fisher-Rao metric (see e.g.\ Section 2.2 in \cite{Gallouet2017} and the references therein), which, for any $\mu_0$, $\mu_1 \in \mathcal{P}_{ac}(\mathbb{R}^d)$ states that 
\[ 
\begin{split}
\operatorname{FR}(\mu_0,\mu_1)=\inf \left\{  \int_0^1 \int_{\mathbb{R}^d} |\nu_s|^2 m_s(dx)ds\,:\, \text{s.t} \,\,\,
\partial_s m_s = \nu_s m_s\,,\quad m_{i}=\mu_i\,,\,\, i=0,1\right\} \,,
\end{split}
\]
where the infimum is taken over all curves $[0,1] \ni t \mapsto (m_t,\nu_t) \in \mathcal{P}_{ac}(\mathbb{R}^d) \times L^2(\mathbb{R}^d;m_t)$ solving $\partial_t m_t = \nu_t m_t$ in the distributional sense, such that $t \mapsto m_t$ is weakly continuous with endpoints $\mu_0$ and $\mu_1$. 
This result tells us that measures in the space $(\mathcal P_{ac}(\mathbb R^d), \operatorname{FR} )$ are transported along curves prescribed by a birth-death (or reaction) equation. 
The main focus of this work is to identify a corresponding Polyak-\L ojasiewicz inequality from which we can deduce the exponential convergence to $m^{\sigma,*}$ of the flow $(m_t)_{t \geq 0}$ described by the birth-death equation
\begin{equation}\label{eq:birthdeath}
\partial_t m_t(x) = - a(m_t,x)
 m_t(x)\,,\,\,\, \quad a(m,x):=
	\frac{\delta F}{\delta m} (m,x) + \frac{\sigma^2}{2} \log \left( \frac{m(x)}{\pi(x)}\right) - \frac{\sigma^2}{2} \operatorname{KL}(m|\pi).
\end{equation}
Note that the map $(m,x)\mapsto a(m,x)$ formally corresponds to $\frac{\delta V^{\sigma}}{\delta m}(m_t,\cdot)$ which may not exist since the KL-divergence is only lower semi-continuous. The map $(m,x)\mapsto a(m,x)$
is a well-defined function under the assumption of flat-differentiability of $F$ (note that $\operatorname{KL}(m|\pi)$ in \eqref{eq:birthdeath} corresponds to the normalizing constant needed in our normalizing convention mentioned above).

To see why the particular form of $(m,x)\mapsto a(m,x)$ in \eqref{eq:birthdeath} is a good choice one needs to show that 
 $t \mapsto V^\sigma(m_t)$ is differentiable so that
 \begin{equation}\label{eq dissipation}
		\begin{split}
		\partial_t &\left( V^{\sigma}(m_t) - V^{\sigma}(m^{\sigma,*})\right)\\
		&=
		\int_{\mathbb{R}^d} \left( \frac{\delta F}{\delta m} (m_t,x) + \frac{\sigma^2}{2}\log \left( \frac{m_t(x)}{\pi(x)}\right) - \frac{\sigma^2}{2}\operatorname{KL}(m_t|\pi) \right) \partial_t m_t(x) dx \\
		 & = - \int_{\mathbb{R}^d} \left| a(m_t,x) \right|^2 m_t(x)dx\,.
		\end{split}
\end{equation}
The Polyak-\L ojasiewicz condition that implies the exponential convergence of $V^{\sigma}(m_t)$ to $V^{\sigma}(m^{\sigma,*})$, requires that there exists a constant $\kappa>0$ such that for any $m^{*}\in \arg \min_m V^{\sigma}(m)$ and any $m\in \mathcal P(\mathbb R^d)$, 
\begin{equation}\label{eq:Lojasiewicz bd}
\frac{1}{ \kappa}\left\| a(m,\cdot)\right\|^2_{L^2(m)} \geq  V^{\sigma}(m) - V^{\sigma}(m^*) \,.
\end{equation}
We call \eqref{eq:Lojasiewicz bd} the flat Polyak-\L ojasiewicz condition, since the function $a(m,x)$ formally corresponds to the flat derivative of $V^{\sigma}$, as explained above. With such an inequality at hand, one immediately sees that 
\begin{equation*}
		\begin{split}
		\partial_t ( V^{\sigma}(m_t) - V(m^{\sigma,*})) 
		\leq  - \kappa (V^{\sigma}(m_t) - V(m^{\sigma,*}) )\,,
		\end{split}
		\end{equation*}
which implies that $V^{\sigma}(m_t) - V^{\sigma}(m^{\sigma,*}) \leq \left( V^{\sigma}(m_0) - V^{\sigma}(m^{\sigma,*}) \right) e^{-\kappa t}$ holds for any $t \geq 0$.

The main contributions of this work are:
\begin{itemize}
\item We establish the existence and uniqueness of the non-linear infinite dimensional birth-death flow \eqref{eq:birthdeath}.
\item We demonstrate that $t \mapsto V^\sigma(m_t)$ is differentiable, which implies that the energy dissipation equality \eqref{eq dissipation} holds. 
\item We show that for a large class of energy functions $V^{\sigma}$, the Polyak-\L ojasiewicz condition \eqref{eq:Lojasiewicz bd} can be verified under relatively mild assumptions.
\end{itemize}

We remark that showing the existence of a solution to \eqref{eq:birthdeath} is non-trivial, since the problem is non-linear and the coefficient $a(m,x)$ contains two terms that are difficult to control: the flat derivative of $F$ and the KL-divergence. Even if one assumes a priori that the former is bounded, it is still unclear how to control the latter. We deal with this problem by introducing a Picard iteration gradient flow approximating \eqref{eq:birthdeath}, and then analysing the symmetrised KL-divergence (rather than just the plain KL-divergence) along that auxiliary gradient flow (see Lemmas \ref{lem:existenceGF} and \ref{thm:Contractivity}). This approach allows us also to obtain bounds on the Radon-Nikodym derivative of $m_t$ with respect to $\pi$ (see Theorem \ref{thm:existence}) that are crucial for proving the differentiability of $t \mapsto V^\sigma(m_t)$ (Theorem \ref{thm:differentiability}) as well as establishing the Polyak-\L ojasiewicz condition \eqref{eq:Lojasiewicz bd} in Theorem \ref{thm:PLI}.

It can be shown that the birth-death flow \eqref{eq:birthdeath} is a limit of a minimising movement scheme, see e.g. \cite{santambrogio2017euclidean},  defined for $\tau > 0$ as 
\begin{equation*}
\mu_{n+1} = \operatorname{argmin}_{\nu \in \mathcal P_{ac}(\mathbb R^d)} \left\{  V^{\sigma}(\nu) + \frac{1}{\tau} \operatorname{KL}(\nu|\mu_n)  \right\} \,.
\end{equation*}
Indeed, recalling that $(m_t,x)\mapsto a(m_t,x)$ defined in \eqref{eq:birthdeath} formally corresponds to $\frac{\delta V^{\sigma}}{\delta m}(m_t,x)$ and using Proposition 2.5 in \cite{HuRenSiskaSzpruch2021} we can easily see that 
\begin{equation}\label{eq:discreteBD}
\frac{\log \mu_{n+1}(x) - \log \mu_n(x)}{\tau} = - a(\mu_{n+1},x) \,.
\end{equation}  
This is an implicit Euler discretisation of \eqref{eq:birthdeath}. Similarly one can consider the mirror descent algorithm, recently studied in \cite{Aubin-Frankowski2022} for the problem of optimization over the space of measures. One can define  
\begin{equation*}
\bar{\mu}_{n+1} = \operatorname{argmin}_{\nu \in P_{ac}(\mathbb R^d)} \left\{  \int_{ R^d} a(\bar{\mu}_n,x) (\nu - \bar{\mu}_n)(dx) +  \frac{1}{\tau} \operatorname{KL}(\nu|\bar{\mu}_n)    \right\} \,.
\end{equation*}
As before, one can show, using Proposition 2.5 in \cite{HuRenSiskaSzpruch2021}, that 
\begin{equation}\label{eq:discreteBD2}
\frac{\log \bar{\mu}_{n+1}(x) - \log \bar{\mu}_n(x)}{\tau} = - a(\bar{\mu}_{n},x) \,.
\end{equation}  
This is an explicit Euler discretisation of \eqref{eq:birthdeath}.  
Note that Theorem 4 in \cite{Aubin-Frankowski2022} shows convergence of their energy function evaluated at $\bar{\mu}_n$ under certain strong convexity assumptions, whereas we work with the (measure space version of) Polyak-\L ojasiewicz inequality. In this context our results provide a natural extension of convergence results for mirror descent algorithms on $\mathbb{R}^d$, which are known to converge under the classical P\L I \eqref{eq:PLIonRd}, see \cite{Radhakrishnan2020}.

The remaining part of the paper is organised as follows. In Section \ref{section:mainResults} we formulate our main results and the assumptions we work with. In Section \ref{section:GeneralPLI} we present a result on the verification of the flat Polyak-\L ojasiewicz inequality \eqref{eq:Lojasiewicz bd} for general energy functions (not necessarily of the form \eqref{eq:Vsigma}) under certain quadratic growth conditions. This section is of independent interest and can be seen as a counterpart of the results that were proved in $\mathbb{R}^d$ in \cite{Karimi2016}, or the results that were proved on the space of measures in \cite{BlanchetBolte2018} for a quadratic growth condition with respect to the $L^2$-Wasserstein distance (while we work with the KL-divergence and the $\chi^2$-divergence). In Section \ref{section:literature} we review the literature and we present a more in-depth discussion on the motivation for studying the gradient flow \eqref{eq:birthdeath}. In Section \ref{section:existence} we prove our main results on the existence of the gradient flow and the differentiability of the energy function. Finally, the Appendix includes some general auxiliary results on comparing different $f$-divergences, adapted from \cite{Dragomir} and a brief overview of the notion of the flat derivative.

\section{Main Results}\label{section:mainResults}

\subsection{Existence of the birth-death flow and its convergence under the flat Polyak-\L ojasiewicz condition}

We work with the energy function $V^{\sigma}: \mathcal{P}(\mathbb{R}^d) \to \mathbb{R}$ given by \eqref{eq:Vsigma}, for some possibly non-linear $F: \mathcal{P}(\mathbb{R}^d) \to \mathbb{R}$ and $\sigma > 0$. We have the following assumptions on $F$.

\begin{assumption}\label{assump:F}
	Suppose $F$ has the first and the second order flat derivatives ($\frac{\delta F}{\delta m} : \mathcal{P}(\mathbb{R}^d) \times \mathbb{R}^d \to \mathbb{R}$ and $\frac{\delta^2 F}{\delta m^2}: \mathcal{P}(\mathbb{R}^d) \times \mathbb{R}^d \times \mathbb{R}^d \to \mathbb{R}$, respectively). Furthermore, suppose that
	\begin{enumerate}
		\item $F$ is convex, i.e., for any $m$, $m' \in \mathcal{P}(\mathbb{R}^d)$ we have
		\begin{equation}\label{eq:convexF}
		F(m) - F(m') \leq \int_{\mathbb{R}^d} \frac{\delta F}{\delta m}(m,x) \left( m - m' \right) (dx) \,.
		\end{equation}
		\item There exists a constant $C > 0$ such that for all $m \in  \mathcal{P}(\mathbb{R}^d)$ and for all $x \in \mathbb{R}^d$ we have
		\begin{equation}\label{eq:assumpF}
		\left| \frac{\delta F}{\delta m} (m,x) \right| \leq C \,.
		\end{equation}
		\item There exists a constant $C_2 > 0$ such that for all $m \in \mathcal{P}(\mathbb{R}^d)$ and for all $x$, $y \in \mathbb{R}^d$ we have
		\begin{equation}\label{eq:assumpF2}
		\left| \frac{\delta^2 F}{\delta m^2} (m,x,y) \right| \leq C_2 \,.
		\end{equation}
	\end{enumerate} 
\end{assumption}

Furthermore, suppose we have absolutely continuous probability measures $\pi$, $m_0 \in \mathcal{P}(\mathbb{R}^d)$ such that $\pi(dx) \propto \exp \left( -\frac{2}{\sigma^2}U(x) \right) dx$ for a potential $U: \mathbb{R}^d \to \mathbb{R}$ and the following conditions are satisfied.

\begin{assumption}\label{assump:m0}
		Suppose $m_0 \in \mathcal{P}(\mathbb{R}^d)$ is absolutely continuous and comparable with $\pi$ in the following sense.
		\begin{enumerate}
		\item There exists a constant $r>0$ such that
		\begin{equation}\label{eq:assumpInf}
		\inf_{x \in \mathbb{R}^d} \frac{m_0(x)}{\pi(x)} \geq r \,.
		\end{equation}
		\item There exists a constant $R > 1$ such that
		\begin{equation}\label{eq:assumpSup}
		\sup_{x \in \mathbb{R}^d} \frac{m_0(x)}{\pi(x)} \leq R \,.
		\end{equation}
	\end{enumerate}
\end{assumption}

Note that here $\pi$ is just a reference measure, and recall that the actual measure of interest (the minimizer of $V^{\sigma}$) is given implicitly by the following equation
\begin{equation*}
m^{\sigma,*}(x) = \frac{1}{Z} \exp{\left( -\frac{2}{\sigma^2} \left( \frac{\delta F}{\delta m} (m^{\sigma,*},x) + U(x) \right)\right)} \,,
\end{equation*}
where $Z$ is the normalizing constant. We immediately observe that, under condition \eqref{eq:assumpF}, conditions \eqref{eq:assumpInf} and \eqref{eq:assumpSup} together are equivalent to assuming that there exist constants $\bar{r}>0$, $\bar{R} > 1$ such that for all $x \in \mathbb{R}^d$,
\begin{equation}\label{eq:ratioStar}
\bar{r} \leq \frac{m_0(x)}{m^{\sigma,*}(x)} \leq \bar{R} \,.
\end{equation}  

As we will explain in more detail in Subsection \ref{section:literature}, Assumption \ref{assump:m0} is a kind of "warm start" condition that says that once we fix the reference measure $\pi$ in \eqref{eq:Vsigma}, the initial measure $m_0$ of our gradient flow should be comparable to $\pi$. 
We have the following result.

\begin{theorem}\label{thm:existence}
	Under Assumption \ref{assump:F} and condition \eqref{eq:assumpSup} from Assumption \ref{assump:m0}, equation \eqref{eq:birthdeath} has a unique solution $(m_t)_{t \geq 0}$. Moreover, for $t\geq 0$,
	\begin{equation}\label{eq:KLmtBound}
\operatorname{KL}(m_t|\pi)  \leq 2 \log R + \frac{4C}{\sigma^2}
\end{equation}
and there exists a constant $R_1 > 1$ such that for all $t \geq 0$,
\begin{equation}\label{eq:ratioBoundAllt}
\sup_{x \in \mathbb{R}^d} \frac{m_t(x)}{\pi(x)} \leq R_1 \,.
\end{equation}
If we additionally assume that condition \eqref{eq:assumpInf} from Assumption \ref{assump:m0} holds, then there exists a constant $r_1 > 0$ such that for all $t \geq 0$,
\begin{equation}\label{eq:ratioBoundAllt2}
\inf_{x \in \mathbb{R}^d} \frac{m_t(x)}{\pi(x)} \geq r_1 \,.
\end{equation}
\end{theorem}

As we explained in the discussion in Section \ref{section:introduction}, the crucial property needed for showing the exponential convergence of $(m_t)_{t \geq 0}$ is the differentiability of the energy function along the gradient flow.

\begin{theorem}\label{thm:differentiability}
		Under Assumption \ref{assump:F} and condition \eqref{eq:assumpSup} from Assumption \ref{assump:m0}, for the unique solution $(m_t)_{t \geq 0}$ to \eqref{eq:birthdeath}, the function $t \mapsto V^{\sigma}(m_t)$ is differentiable and
		\begin{equation}\label{eq:diff}
		\partial_t V^{\sigma}(m_t) = - \int_{\mathbb{R}^d} \left| 	\frac{\delta F}{\delta m} (m_t,x) + \frac{\sigma^2}{2} \log \left( \frac{m_t(x)}{\pi(x)}\right) - \frac{\sigma^2}{2} \operatorname{KL}(m_t|\pi) \right|^2 m_t(x) dx \,. 
		\end{equation}
\end{theorem}

Note that inequalities \eqref{eq:ratioBoundAllt} and \eqref{eq:ratioBoundAllt2} obtained in Theorem \ref{thm:existence} imply that there exist constants $\bar{r}_1 > 0$ and $\bar{R}_1 > 1$ are such that for all $t \geq 0$ and all $x \in \mathbb{R}^d$,
\begin{equation*}
\bar{r}_1 \leq \frac{m_t(x)}{m^{\sigma,*}(x)} \leq \bar{R}_1 
\end{equation*} 
(similarly to how \eqref{eq:assumpInf} and \eqref{eq:assumpSup} imply \eqref{eq:ratioStar}). This property will be crucial in the proof of the following Polyak-\L ojasiewicz inequality.
\begin{theorem}\label{thm:PLI}
	Under Assumptions \ref{assump:F} and \ref{assump:m0}, the flow $(m_t)_{t \geq 0}$ solving \eqref{eq:birthdeath} satisfies 
	\begin{equation}\label{eq:flatLojasiewicz2}
	V^{\sigma}(m_t) - V^{\sigma}(m^{\sigma,*}) \leq \frac{4 \bar{R}_1}{\sigma^2 \bar{r}_1} \left\| a(m_t, \cdot) \right\|^2_{L^2(m_t)} 
	\end{equation}
	for all $t \geq 0$.
\end{theorem}

Using Theorem \ref{thm:PLI}, based on the discussion in Section \ref{section:introduction}, we have the following result.

\begin{corollary}
		Under Assumptions \ref{assump:F} and \ref{assump:m0}, the flow $(m_t)_{t \geq 0}$ solving \eqref{eq:birthdeath} satisfies 
		\begin{equation*}
		V^{\sigma}(m_t) - V^{\sigma}(m^{\sigma,*}) \leq \left( V^{\sigma}(m_0) - V^{\sigma}(m^{\sigma,*}) \right) e^{-\kappa t} \,,
		\end{equation*}
		for all $t \geq 0$, where $\kappa = \sigma^2 \bar{r} / 4 \bar{R}$. 
\end{corollary}

The proofs of all the results formulated above are postponed to Section \ref{section:existence}.

In Subsection \ref{section:GeneralPLI} we will explain how to deduce the Polyak-\L ojasiewicz inequality \eqref{eq:flatLojasiewicz2} for a general class of energy functions that satisfy a certain growth condition with respect to the KL-divergence. We will now formulate a lemma where we verify that growth condition for the energy function $V^{\sigma}$ that we used in this subsection (given by \eqref{eq:Vsigma}).

\begin{lemma}\label{lem:quadraticGrowth}
	For $V^{\sigma}$ given by \eqref{eq:Vsigma}, if $F$ is convex, then $V^{\sigma}$ satisfies the quadratic growth condition 
	\begin{equation*}
	V^{\sigma}(m) - V^{\sigma}(m^{\sigma,*}) \geq \frac{\sigma^2}{2} \operatorname{KL}(m | m^{\sigma,*}) 
	\end{equation*}
	for any $m \in \mathcal{P}(\mathbb{R}^d)$.
		\begin{proof}
		The proof is a straightforward extension of the proof of Proposition 1 in \cite{Nitanda2022}, where this was shown for $V = F + H$, where $H$ is the negative entropy. 
		By convexity of $F$, for any probability measures $m$, $m' \in \mathcal{P}(\mathbb{R}^d)$ we get
		\begin{equation*}
		\begin{split}
		&V^{\sigma}(m') = F(m') + \frac{\sigma^2}{2} \operatorname{KL} (m' | \pi)\\
		& \geq F(m) + \int_{\mathbb{R}^d} \frac{\delta F}{\delta m} (m,x) (m'-m)(dx) + \frac{\sigma^2}{2} \operatorname{KL} (m' | \pi)\\
		&= F(m) + \int_{\mathbb{R}^d} \left(\frac{\delta F}{\delta m} (m,x) + \frac{\sigma^2}{2} \log \frac{m(x)}{\pi(x)} - \frac{\sigma^2}{2} \log \frac{m(x)}{\pi(x)} \right) (m'-m)(dx) + \frac{\sigma^2}{2} \operatorname{KL} (m' | \pi)\\
		& = F(m) + \int_{\mathbb{R}^d} a(m,x) (m'-m)(dx)  -\int_{\mathbb{R}^d} \frac{\sigma^2}{2} \log \frac{m(x)}{\pi(x)} (m'-m)(dx)+ \frac{\sigma^2}{2} \operatorname{KL} (m' | \pi)\\
		& = F(m) + \int_{\mathbb{R}^d} a(m,x) (m'-m)(dx) + \frac{\sigma^2}{2} \operatorname{KL}(m'|m) + \frac{\sigma^2}{2} \operatorname{KL}(m|\pi)\\
		& \geq V^{\sigma}(m) + \int_{\mathbb{R}^d} a(m,x) (m'-m)(dx) + \frac{\sigma^2}{2} \operatorname{KL}(m'|m) \,.
		\end{split}
		\end{equation*}
		Taking $m=m^{\sigma,*}$ in the above calculation finishes the proof, since $a(m^{\sigma,*},\cdot)$ is constant by Proposition 2.5 in \cite{HuRenSiskaSzpruch2021}.
	\end{proof}	
\end{lemma}	

Note that we call the growth condition in Lemma \ref{lem:quadraticGrowth} quadratic, since the KL-divergence corresponds to the square of a distance on the space of measures (compare this to condition (2) for $\theta = 1/2$ in \cite{BlanchetBolte2018}, which considered a similar growth condition with the $L^2$-Wasserstein distance, and see the discussion below our Remark \ref{remark:kappa} for more details).

 \subsection{Verification of the flat Polyak-\L ojasiewicz condition in a general setting}\label{section:GeneralPLI}

In this subsection we adapt the proof of Theorem 2 in \cite{Karimi2016} to the setting of the space of measures. In \cite{Karimi2016} it was shown how the classical Polyak-\L ojasiewicz inequality \eqref{eq:PLIonRd} for functions on $\mathbb{R}^d$ can be inferred from a certain type of a quadratic growth condition. Here we will work with functions on $\mathcal{P}(\mathbb{R}^d)$ and we will carry out a similar argument, based on certain quadratic growth conditions expressed in terms of either the KL-divergence or the $\chi^2$-divergence, where the latter is defined 
for any $m \in \mathcal{P}(\mathbb{R}^d)$ by
\begin{equation*}
\chi^2(m|\pi) =\begin{cases}
\int_{\mathbb{R}^d} \left( \frac{m(x)}{\pi(x)} - 1\right)^2 \pi(x) dx \, &\text{$m$ absolutely continuous with respect to $\pi$,}  \\
\infty  &\text{otherwise.}
\end{cases}
\end{equation*} 
This result can be interpreted as an analogue of Theorem 1 in \cite{BlanchetBolte2018}, which showed that a certain type of the \L ojasiewicz inequality can be inferred from a quadratic growth condition with respect to the $L^2$-Wasserstein distance. We will present our reasoning in a series of lemmas.

\begin{lemma}\label{lemma:flatCauchySchwarz}
	Suppose that $G: \mathcal{P}(\mathbb{R}^d) \to \mathbb{R}$ has the first order flat derivative and that $G$ is convex (cf.\ \eqref{eq:convexF}). 
	Then for any absolutely continuous probability measures $m$, $m' \in \mathcal{P}(\mathbb{R}^d)$,
	\begin{equation*}
	\begin{split}
	G(m) - G(m') &\leq \left( \int_{\mathbb{R}^d} \left| \frac{\delta G}{\delta m}(m,x) \right|^2 m(x) dx \right)^{1/2} \left( \int_{\mathbb{R}^d} \left( \frac{m'(x)}{m(x)} - 1 \right)^2 m(x) dx \right)^{1/2} \\
	&= \left\| \frac{\delta G}{\delta m} (m,\cdot) \right\|_{L^2(m)} \cdot \chi^2(m'|m)^{1/2} \,.
	\end{split}
	\end{equation*}
	\begin{proof}
		Since $\int_{\mathbb{R}^d} \frac{\delta G}{\delta m}(m,x)  m(x) dx = 0$ by convention, from the convexity condition \eqref{eq:convexF} we get
		\begin{equation*}
		G(m) - G(m') \leq 
		- \int_{\mathbb{R}^d} \frac{\delta G}{\delta m}(m,x)  \left( \frac{m'(x)}{m(x)} - 1 \right) m(x)dx \,.
		\end{equation*}
		A simple application of the Cauchy-Schwarz inequality in $L^2(m)$ proves the desired assertion.
	\end{proof}
\end{lemma}

Next we need a lemma that allows us to compare the $\chi^2$-divergence and the KL-divergence, between two absolutely continuous measures, such that the ratio of their densities is bounded from above and below.

\begin{lemma}\label{lemma:Dragomir}
	Suppose we have absolutely continuous $m$, $m' \in \mathcal{P}(\mathbb{R}^d)$ such that there exist constants $r$, $R > 0$ such that for any $x \in \mathbb{R}^d$ we have
	\begin{equation*}
	r \leq \frac{m(x)}{m'(x)} \leq R \,.
	\end{equation*}
	Then we have
	\begin{equation}\label{eq:Dragomir}
	\operatorname{KL}(m'|m) \leq \frac{1}{r} \operatorname{KL}(m|m') \qquad \text{ and } \qquad \chi^2(m|m') \leq 2R \operatorname{KL}(m|m') \,.
	\end{equation}
	\begin{proof}
		The proof can be adapted from the proofs of Proposition 1 and Proposition 2 in \cite{Dragomir}, which covered the case of discrete probability measures. For completeness, we include the proof in Section \ref{section:appendix1}. 
	\end{proof}
\end{lemma}

Based on the above lemmas, we can show the following result.

\begin{theorem}\label{thm:mainLojasiewicz}
	Suppose that $G: \mathcal{P}(\mathbb{R}^d) \to \mathbb{R}$ has the first order flat derivative and that $G$ is convex. Suppose further that $G$ is minimized by an absolutely continuous measure $m^*$ and that there exists a constant $\lambda > 0$ such that for any $m' \in \mathcal{P}(\mathbb{R}^d)$,
	\begin{equation}\label{eq:quadraticGrowthKL}
	G(m') - G(m^*) \geq \lambda \operatorname{KL}(m'|m^*) \,.
	\end{equation}
	Moreover, suppose that we have an absolutely continuous measure $m \in \mathcal{P}(\mathbb{R}^d)$ such that there exist constants $r$, $R > 0$ such that for any $x \in \mathbb{R}^d$ we have
	\begin{equation}\label{eq:ratio}
	r \leq \frac{m(x)}{m^*(x)} \leq R \,.
	\end{equation}
	Then 
	\begin{equation}\label{eq:PLIG}
	G(m) - G(m^*) \leq \frac{2R}{\lambda r} \left\| \frac{\delta G}{\delta m} (m,\cdot) \right\|^2_{L^2(m)} \,.
	\end{equation}
	\begin{proof}
		We follow the argument from the proof of Theorem 1 in \cite{BlanchetBolte2018}. 
		Since $G$ is assumed to be convex, from Lemma \ref{lemma:flatCauchySchwarz} we get  
		\begin{equation}\label{eq:newLojAux}
		G(m) - G(m^*) \leq \left\| \frac{\delta G}{\delta m} (m,\cdot) \right\|_{L^2(m)} \cdot \chi^2(m^*|m)^{1/2} \,.
		\end{equation}
		However, due to Lemma \ref{lemma:Dragomir}, we have
		\begin{equation*}
		\chi^2(m^*|m) \leq 2R \operatorname{KL}(m^*|m) \leq \frac{2R}{r} \operatorname{KL}(m|m^*) \,,
		\end{equation*}
		which, together with \eqref{eq:newLojAux} and $G(m) - G(m^*) \geq \lambda \operatorname{KL}(m|m^*)$ leads to
		\begin{equation*}
		\operatorname{KL}(m|m^*)^{1/2} \leq \frac{1}{\lambda} \left( \frac{2R}{r} \right)^{1/2} \left\| \frac{\delta G}{\delta m} (m,\cdot) \right\|_{L^2(m)} \,.
		\end{equation*}
		In particular,
		\begin{equation}\label{eq:newLojAux2}
		\chi^2(m^*|m)^{1/2} \leq \frac{2R}{\lambda r} \left\| \frac{\delta G}{\delta m} (m,\cdot) \right\|_{L^2(m)} \,.
		\end{equation}
		Plugging \eqref{eq:newLojAux2} into the right hand side of \eqref{eq:newLojAux}, we obtain
		\begin{equation*}
		G(m) - G(m^*) \leq \frac{2R}{\lambda r} \left\| \frac{\delta G}{\delta m} (m,\cdot) \right\|_{L^2(m)}^2 \,.
		\end{equation*}
	\end{proof}
\end{theorem}

\begin{remark}\label{remark:kappa}
	Under the assumptions of Theorem \ref{thm:mainLojasiewicz}, we obtain the flat Polyak-\L ojasiewicz condition of the type \eqref{eq:Lojasiewicz bd} with the constant 
	\begin{equation}\label{eq:kappa}
	\kappa = \left(  \frac{2R}{\lambda r} \right)^{-1} \,.
	\end{equation}
	In what follows, we will prove that the flow $(m_t)_{t \geq 0}$ given by \eqref{eq:birthdeath} is such that $\bar{r}_1 \leq \frac{m_t(x)}{m^{\sigma,*}(x)} \leq \bar{R}_1$ with some constants $\bar{r}_1 > 0$, $\bar{R}_1 > 1$, for all $t >0$ and $x \in \mathbb{R}^d$, which will allow us to show \eqref{eq:PLIG} with $G$ on the left hand side replaced by $V^{\sigma}$, and $\frac{\delta G}{\delta m}(m,x)$ on the right hand side replaced by $a(m,x)$ given by \eqref{eq:birthdeath}. This will be the basis of the proof of our main results in Section \ref{section:existence} and will provide us with an exponential convergence rate of $V^{\sigma}(m_t)$ to $V^{\sigma}(m^{\sigma,*})$. 
 We can easily observe that the convergence rate $\kappa$ given by \eqref{eq:kappa} degenerates to zero when $\lambda \to 0$ or $r \to 0$ or $R \to \infty$.
\end{remark}

Condition \eqref{eq:quadraticGrowthKL} corresponds to the classical quadratic growth condition for functions $f: \mathbb{R}^d \to \mathbb{R}$ that can be used (see Theorem 2 in \cite{Karimi2016}) to prove the classical Polyak-\L ojasiewicz inequality \eqref{eq:PLIonRd} under the additional assumption of convexity of $f$ (but not necessarily strong convexity). More precisely, the quadratic growth condition in $\mathbb{R}^d$ states that
\begin{equation*}
f(x) - \min_{y\in \mathbb R^d} f(y) \geq \frac{\mu}{2} \| x - x_p \|^2 \,,
\end{equation*}
where $x_p \in \arg \min_{x \in \mathbb{R}^d} f(x)$. Specifying an analogous condition for functions on the space of measures is non-straightforward, as there are multiple choices of the notion of the distance. Blanchet and Bolte in \cite{BlanchetBolte2018} proved that a certain type of a \L ojasiewicz inequality can be implied by a condition such as \eqref{eq:quadraticGrowthKL} but with the $L^2$-Wasserstein distance instead of the KL-divergence, see formula (2) and Theorem 1 in \cite{BlanchetBolte2018}. Based on the proof of our Theorem \ref{thm:mainLojasiewicz}, it is clear that we can also consider a quadratic growth condition with respect to the $\chi^2$-divergence with reversed arguments, i.e., we have the following result.

\begin{corollary}
	Suppose that $G: \mathcal{P}(\mathbb{R}^d) \to \mathbb{R}$ has the first order flat derivative and that $G$ is convex. Suppose further that $G$ is minimized by an absolutely continuous measure $m^*$ and that there exists a constant $\lambda > 0$ such that for any $m' \in \mathcal{P}(\mathbb{R}^d)$,
	\begin{equation}\label{eq:quadraticGrowthChi2}
	G(m') - G(m^*) \geq \lambda \chi^2(m^*|m') \,.
	\end{equation}
Then for any $m \in \mathcal{P}(\mathbb{R}^d)$ we have the flat Polyak-\L ojasiewicz condition 
		\begin{equation}\label{eq:PLIG2}
G(m) - G(m^*) \leq \frac{1}{\lambda} \left\| \frac{\delta G}{\delta m} (m,\cdot) \right\|_{L^2(m)}^2 \,.
\end{equation}
	\begin{proof}
Using \eqref{eq:newLojAux} and \eqref{eq:quadraticGrowthChi2}, one immediately obtains  
		\begin{equation*}
		\chi^2(m^*|m)^{1/2} \leq \frac{1}{\lambda} \left\| \frac{\delta G}{\delta m} (m,\cdot) \right\|_{L^2(m)} \,,
		\end{equation*}
		which can be plugged back into \eqref{eq:newLojAux} to obtain \eqref{eq:PLIG2}.
	\end{proof}	
\end{corollary}

The quadratic growth condition with respect to the KL-divergence \eqref{eq:quadraticGrowthKL} seems more natural than the one with respect to the $\chi^2$-divergence \eqref{eq:quadraticGrowthChi2} (note that the former is verified in Lemma \ref{lem:quadraticGrowth} for a large class of energy functions given by \eqref{eq:Vsigma}). It is clear based on Lemma \ref{lemma:Dragomir} that \eqref{eq:quadraticGrowthKL} implies \eqref{eq:quadraticGrowthChi2}, but we are presently unaware of any examples of energy functions that would satisfy \eqref{eq:quadraticGrowthChi2} but not \eqref{eq:quadraticGrowthKL}.

\subsection{Literature review, connection to the Wasserstein-Fisher-Rao gradient flow and further research }\label{section:literature}

In order to present our results in a broader context, let us first discuss a different type of gradient flows and associated \L ojasiewicz-type inequalities. We will also provide two heuristic examples in order to build a better intuition for our approach.

\subsubsection{Wasserstein Gradient Flow}

The dynamic representation of the $L^2$-Wasserstein metric $\mathcal W_2$ due to Benamou and Brenier \cite{BenamouBrenier2000, Villani} states that for any $\mu_0$, $\mu_1 \in \mathcal P_2(\mathbb R^d)$, 
\begin{equation}\label{eq:Benamou}
\mathcal W_2(\mu_0,\mu_1) 
 =\inf \left\{  \int_0^1 \int_{\mathbb{R}^d} |\nu_s|^2  m_s(dx)ds :\, \text{s.t} \,\,
\partial_s m_s + \div  (\nu_s m_s) =0\,,\, m_{i}=\mu_i\,,\,\, i=0,1 \right\},
\end{equation}
where the infimum is taken over all curves $[0,1] \ni t \mapsto (m_t,\nu_t) \in \mathcal P_2(\mathbb R^d) \times L^2(\mathbb{R}^d;m_t)$ solving $\partial_t m_t + \div  (\nu_t m_t) =0$ in the distributional sense, such that $t \mapsto m_t$ is weakly continuous with endpoints $\mu_0$ and $\mu_1$. 
This result tells us that measures in the space $(\mathcal P_2(\mathbb R^d), \mathcal W_2 )$ of probability measures with finite second moments are transported along curves described by the forward-Kolmogorov PDE.  

One can show \cite{HuRenSiskaSzpruch2021} that  $V^{\sigma}(m_t)$ is decreasing along the gradient flow $(m_t)_{t \geq 0}$ satisfying 
\begin{equation}\label{eq:gradientFlow}
\partial_t m_t     = \div  \left( \nabla a(m_t,\cdot)m_t\right)
\,,\,\,\,\,\,  a(m,x):=
	\frac{\delta F}{\delta m} (m,x) + \frac{\sigma^2}{2} \log \left( \frac{m(x)}{\pi(x)}\right) - \frac{\sigma^2}{2} \operatorname{KL}(m|\pi).
\end{equation}
Note that this flow corresponds to the mean-field Langevin equation (see e.g.\ (1.4) and (1.5) in \cite{HuRenSiskaSzpruch2021}), and in particular becomes the classical overdamped Langevin equation when $F=0$. 
Indeed, if we can show that $t \mapsto V^{\sigma}(m_t)$ is differentiable (see e.g.\ \cite[Theorem  2.9]{HuRenSiskaSzpruch2021}), we obtain 
\begin{equation}\label{eq:gfCalc}
		\begin{split}
		\partial_t V^{\sigma}(m_t) &=
		\int_{\mathbb{R}^d} a(m_t,x) \partial_t m_t(x) dx = \int_{\mathbb{R}^d} a(m_t,x) \div  \left( \nabla a(m_t,x)m_t(x)\right) dx  \\
		&= - \int_{\mathbb{R}^d} \left| \nabla a(m_t,x)\right|^2 m_t(dx)  \,.
		\end{split}
		\end{equation}	
In the case when $F$ is convex, and hence $V^{\sigma}$ is strictly convex, $V^{\sigma}(m_t) \rightarrow V^{\sigma}(m^{\sigma,*})$, see \cite{HuRenSiskaSzpruch2021}. More recently, \cite{Nitanda2022} and \cite{Chizat2022} under additional structural assumptions proved that this convergence is exponential. 

In this setting, the Polyak-\L ojasiewicz condition that implies the exponential convergence $V^{\sigma}(m_t) \to V^{\sigma}(m^{\sigma,*})$, requires that there exists a constant $\kappa>0$ such that for any $m^{*}\in \arg \min_{m \in \mathcal{P}(\mathbb{R}^d)} V^{\sigma}(m)$ and any $t\geq 0$,
\begin{equation}\label{eq:Lojasiewicz}
\frac{1}{ \kappa}\left\| \nabla a (m_t,\cdot) \right\|^2_{L^2(m_t)} \geq V^{\sigma}(m_t) - V^{\sigma}(m^{\sigma,*}) \,.
\end{equation}
With such an inequality at hand, one immediately sees that 
\begin{equation*}
		\begin{split}
		\partial_t ( V^{\sigma}(m_t) - V^{\sigma}(m^{\sigma,*})) 
		&= - \int_{\mathbb{R}^d} \left| \nabla a(m_t,x)\right|^2 m_t(dx) 
		\leq  - \kappa (V^{\sigma}(m_t) - V^{\sigma}(m^{\sigma,*}) )\,,
		\end{split}
		\end{equation*}
and the exponential convergence follows due to the Gronwall lemma. 
\begin{example}\label{exampleFzero}
Let $F=0$ in \eqref{eq:Vsigma}. In this case the minimizing probability measure $m^{\sigma,*} = \arg\min_m V^{\sigma}(m) = \pi$. Then, assuming that we can show that $t \mapsto \operatorname{KL}(m_t|\pi)$ is differentiable, we have
\begin{equation*}
	\partial_t ( V^{\sigma}(m_t) - V^{\sigma}(m^{\sigma,*})) =  \frac{\sigma^2}{2}  \, \partial_t  \operatorname{KL}(m_t|\pi) = - \frac{\sigma^4}{4}  \int_{\mathbb{R}^d} \left| \nabla \, \log \frac{m_t(x)}{\pi(x)} \right|^2 m_t(dx)\,.
\end{equation*}
In this case the Polyak-\L ojasiewicz inequality is just the well-known log-Sobolev inequality
\begin{equation}\label{logSobolev}
\frac{1}{ \kappa}\int_{\mathbb{R}^d} \left| \nabla \, \log \frac{m_t(x)}{\pi(x)} \right|^2 m_t(dx) \geq \operatorname{KL}(m_t|\pi) \,.
\end{equation}
\end{example}

\begin{example}
	Let us consider an example with a different type of energy function. 
	Consider $V^{\sigma}(m) := \chi^2(m|\pi) = \int_{\mathbb{R}^d} \left( \frac{m(x)}{\pi(x)} - 1 \right)^2 \pi(x)dx$ for probability measures $m \in \mathcal{P}(\mathbb{R}^d)$ absolutely continuous with respect to $\pi$, and denote 
	 \[
 \bar{a}(m,x)
	:= 2 \left( \frac{m(x)}{\pi(x)} - 1 \right) - 2 \chi^2(m|\pi) \,,
	\] 
	which formally corresponds to the flat derivative of the $\chi^2$-divergence. Then $V^{\sigma}(m_t)$ is decreasing along the gradient flow $(m_t)_{t \geq 0}$ satisfying 
	\begin{equation}\label{eq:gradientFlow2}
	\partial_t m_t     = \div  \left( \nabla \bar{a}(m_t,\cdot)\pi\right)
	\,,
	\end{equation}
	i.e., similarly as in \eqref{eq:gfCalc}, assuming $t \mapsto \chi^2(m_t|\pi)$ is differentiable, we have
	\begin{equation*}
	\partial_t V^{\sigma}(m_t) = - \int_{\mathbb{R}^d} \left| \nabla \bar{a}(m_t,x)\right|^2 \pi(dx) \,.
	\end{equation*}
	Here the Polyak-\L ojasiewicz inequality becomes the Poincar\'{e} inequality
	\begin{equation*}
	\frac{1}{ \kappa}\int_{\mathbb{R}^d} \left| \nabla \left(  \frac{m_t(x)}{\pi(x)} \right) \right|^2 \pi(dx) \geq \chi^2(m_t|\pi) \,.
	\end{equation*}
	Note that this corresponds to \eqref{eq:Lojasiewicz} with the $L^2(\pi)$ norm instead of $L^2(m_t)$, since we used a different gradient flow (compare \eqref{eq:gradientFlow2} to \eqref{eq:gradientFlow}).
\end{example}

\subsubsection{Wasserstein-Fisher-Rao Gradient Flow} 

A natural idea is to combine the Wasserstein \eqref{eq:gradientFlow} and the Fisher-Rao \eqref{eq:birthdeath} gradient flows which in our setting leads to
\begin{equation}\label{eq:WFRflow}
\partial_t m_t =  \div  \left( \nabla a(m_t,\cdot)m_t\right) - a(m_t,x) m_t
\,.
\end{equation}
Flows of this type have been the subject of intensive research over the last few years \cite{Liero2018, Gallouet2017, LuLuNolen2019, vandenEijnden2019}. If we can show the existence of such a flow, and the differentiability of $t \mapsto V^{\sigma}(m_t)$, one can then check that 
\begin{equation*}
\begin{split}
\partial_t \left( V^{\sigma}(m_t) - V^{\sigma}(m^{\sigma,*})\right)
&= - \left\| \nabla a(m_t,\cdot) \right\|^2_{L^2(m_t)} - \left\| a (m_t,\cdot) \right\|^2_{L^2(m_t)} \,.
\end{split}
\end{equation*}
If the corresponding Polyak-\L ojasiewicz conditions \eqref{eq:Lojasiewicz} and \eqref{eq:flatLojasiewicz2} are satisfied, then the right hand side is bounded by $- \left( \sigma^2 \bar{r}_1/(4\bar{R}_1) + \kappa \right) \left( V^{\sigma}(m_t) - V^{\sigma}(m^{\sigma,*}) \right)$ and we easily obtain the exponential convergence $V^{\sigma}(m_t) - V^{\sigma}(m^{\sigma,*}) \leq \left( V^{\sigma}(m_0) - V^{\sigma}(m^{\sigma,*})\right)e^{-\kappa_1 t}$, where $\kappa_1 = \sigma^2 \bar{r}_1/(4\bar{R}_1) + \kappa$.
This shows that both the Langevin part and the birth-death part can independently contribute to the convergence of $V^{\sigma}(m_t)$, if the right corresponding conditions \eqref{eq:Lojasiewicz} or \eqref{eq:flatLojasiewicz2} are satisfied. However, the issues of the existence of $(m_t)_{t \geq 0}$, the differentiability of $t \mapsto V^{\sigma}(m_t)$ and the verification of \eqref{eq:Lojasiewicz} in general settings are all non-trivial and will be studied in our future research, together with the issue of particle system approximation of \eqref{eq:WFRflow}, see also the last paragraph of this section.

We note that \cite{vandenEijnden2019} studied convergence of flows similar to \eqref{eq:WFRflow}. However, they covered energy functions of a very specific form (see (11) in \cite{vandenEijnden2019}) and without regularisation by the KL-divergence. Moreover, \cite{vandenEijnden2019}  obtained an asymptotic polynomial convergence rate in their main result (Theorem 4.6) and they did not address  some important technical issues such as the question of the existence of the gradient flow and the differentiability of $t \mapsto V^{\sigma}(m_t)$.

On the other hand, \cite{LuLuNolen2019} studied \eqref{eq:WFRflow} corresponding to the linear case ($F=0$) of our Example \ref{exampleFzero} and obtained an exponential rate of convergence to $\pi$, measured in the KL-divergence (see Theorem 3.3 therein). Interestingly, even though the authors of \cite{LuLuNolen2019} did not explicitly make a connection to the Polyak-\L ojasiewicz inequalities, their proof is in fact based on showing a special case of condition \eqref{eq:Lojasiewicz bd} as specified above (see their inequality $(2-2\delta)H_1(f) \leq H_2(f)$ in the proof of Theorem 3.3, integrate it with respect to $\rho_t$ and note that our $m_t$ corresponds to their $\rho_t$). This Polyak-\L ojasiewicz inequality is verified in \cite{LuLuNolen2019} under a positive lower bound on the ratio of densities $\inf_{x \in \mathbb{R}^d} \frac{\rho_t(x)}{\pi(x)}$ that is required to hold for all sufficiently large $t$, see (B.3) in \cite{LuLuNolen2019}. Then they use an argument based on the maximum principle (which is possible due to the Langevin component of their dynamics) to show that this condition in fact only has to hold at an initial time $t_0$. As a consequence, they conclude that compared to the classical result on the exponential convergence of the Langevin dynamics to $\pi$ under the log-Sobolev inequality, by adding the birth-death component to the dynamics they can get rid of the log-Sobolev assumption and replace it by a "warm start" condition $\inf_{x \in \mathbb{R}^d} \frac{\rho_{t_0}(x)}{\pi(x)} \geq c$ for some $c >0$. However, in \cite{LuLuNolen2019} the Langevin part of the dynamics is only applied to make the use of the maximum principle possible, and does not directly contribute to the convergence rate. Moreover, similarly as in \cite{vandenEijnden2019}, the question of the existence of the gradient flow and the differentiability of $t \mapsto V^{\sigma}(m_t)$ were not addressed in \cite{LuLuNolen2019}.

In this paper we study a more general setting than \cite{LuLuNolen2019}, including non-linear functions $F$ in the energy function $V^{\sigma}$ in \eqref{eq:Vsigma}, and we rigorously prove the existence of the corresponding birth-death gradient flow $(m_t)_{t \geq 0}$, as well as the differentiability of $t \mapsto V^{\sigma}(m_t)$. We also verify the flat Polyak-\L ojasiewicz inequality \eqref{eq:Lojasiewicz bd} and thus establish the exponential rate of convergence of  $V^{\sigma}(m_t)$ to $V^{\sigma}(m^{\sigma,*})$. Our condition guaranteeing that \eqref{eq:Lojasiewicz bd} holds (Assumption \ref{assump:m0}) resembles the warm start condition from \cite{LuLuNolen2019}, however, in order to show that it propagates from $t=0$ to all $t>0$, we do not use the Langevin component of the dynamics and hence we work with a "pure" birth-death dynamics (the Fisher-Rao gradient flow).

Other recent papers studying the mean-field optimization problem specified by \eqref{eq:Vsigma}, such as \cite{Nitanda2022} and \cite{Chizat2022}, focused on the Wasserstein gradient flow \eqref{eq:gradientFlow}. Both \cite{Nitanda2022} and \cite{Chizat2022} proved the exponential convergence rate of $V^{\sigma}(m_t)$ to $V^{\sigma}(m^{\sigma,*})$ under the assumption of the log-Sobolev inequality for a class of proximal Gibbs measures related to $m^{\sigma,*}$. Compared to \cite{Nitanda2022,Chizat2022}, working with the Fisher-Rao gradient flow allows us to get rid of that assumption, at the cost of introducing the additional "warm start" conditions in Assumption \ref{assump:m0}.

With all that said, we would like to point out that from the point of view of practical algorithms (that will be the subject of our future work), combining the birth-death dynamics with the Langevin dynamics seems advisable. The Wasserstein-Fisher-Rao gradient flow \eqref{eq:WFRflow} can be seen as the mean-field limit of an interacting particle system that can be used as a basis of practically implementable algorithms (as studied in Sections 6 in \cite{LuLuNolen2019} and \cite{vandenEijnden2019}). The support of the birth-death flow does not change in time and hence, intuitively, if we do not include the diffusion component in our dynamics and we initialize it with the empirical measure of a set of particles, the dynamics will just keep re-arranging the mass between the particles but will not change their positions. Hence the convergence of such dynamics should be expected to be worse than the convergence of a particle system utilizing both the Langevin and the birth-death components. This issue is not apparent in the analysis of the mean-field limit process in the present paper (as our results use a "warm start" assumption on the initial condition), but we will investigate it in detail in our future work on the particle system approximations and the corresponding algorithms. From the practical point of view, the main message of this paper is that the birth-death component of such algorithms can be defined in terms of the function $a$ given by \eqref{eq:birthdeath}, which corresponds to the flat derivative of the energy function $V^{\sigma}$, but the focus here is on the  theoretical analysis of the gradient flow rather than applications.

\section{Existence of the gradient flow and other proofs}\label{section:existence}

In order to prove the existence of a solution $(m_t)_{t \geq 0}$ to
\begin{equation}\label{eq:bd1}
\partial_t m_t(x) = - \left( \frac{\delta F}{\delta m} (m_t,x) + \frac{\sigma^2}{2}\log \left( \frac{m_t(x)}{\pi(x)}\right) - \frac{\sigma^2}{2}\operatorname{KL}(m_t|\pi) \right) m_t(x) \,,
\end{equation}
we first notice that \eqref{eq:bd1} is equivalent to
\begin{equation}\label{eq:bd2}
\partial_t \log m_t(x) = - \left( \frac{\delta F}{\delta m}(m_t,x) + \frac{\sigma^2}{2} \log \left( \frac{m_t(x)}{\pi(x)} \right) - \frac{\sigma^2}{2} \operatorname{KL}(m_t|\pi)\right) \,. 
\end{equation}
By Duhamel's formula, \eqref{eq:bd2} is equivalent to
\begin{equation*}
\log m_t(x) = e^{-\frac{\sigma^2}{2}t} \log m_0(x) - \int_0^t \frac{\sigma^2}{2} e^{-\frac{\sigma^2}{2}(t-s)} \left( \frac{2}{\sigma^2} \frac{\delta F}{\delta m}(m_s,x) - \log \pi(x) - \operatorname{KL}(m_s|\pi) \right) ds \,.
\end{equation*}
Based on this formula, we will define a Picard iteration scheme. To this end, let us first fix $T > 0$ and choose a flow of probability measures $(m_t^{(0)})_{t \in [0,T]}$ such that
\begin{equation}\label{eq:KLcon1}
\int_0^T \operatorname{KL}(m_s^{(0)}|\pi) ds < \infty \,.
\end{equation}
For each $n \geq 1$, we want to fix $m_0^{(n)} = m_0^{(0)} = m_0$ (with $m_0$ satisfying condition \eqref{eq:assumpSup} from Assumption \ref{assump:m0}) and define $(m_t^{(n)})_{t \in [0,T]}$ by
\begin{equation}\label{eq:Picard}
\begin{split}
\log m_t^{(n)}(x) &= e^{-\frac{\sigma^2}{2}t} \log m_0(x) \\
&- \int_0^t \frac{\sigma^2}{2} e^{-\frac{\sigma^2}{2}(t-s)} \left( \frac{2}{\sigma^2} \frac{\delta F}{\delta m}(m_s^{(n-1)},x) - \log \pi(x) - \operatorname{KL}(m_s^{(n-1)}|\pi) \right) ds \,.
\end{split}
\end{equation}
We have the following result.

\begin{lemma}\label{lem:existenceGF}
	The sequence of flows $\left( (m_t^{(n)})_{t \in [0,T]} \right)_{n=0}^{\infty}$ given by \eqref{eq:Picard} is well-defined and such that for all $n \geq 1$ and all $t \in [0,T]$ we have
	\begin{equation*}
		\operatorname{KL}(m_t^{(n)}|\pi) \leq 2 \log R + \frac{4}{\sigma^2}C \,.
	\end{equation*}
	\begin{proof}
		Consider $n=1$. By \eqref{eq:assumpF} and \eqref{eq:KLcon1}, the integral on the right hand side of \eqref{eq:Picard} is finite, and hence $(m_t^{(1)})_{t \in [0,T]}$ is well-defined. Note that due to \eqref{eq:assumpF}, the only potential issue with the definition of $(m_t^{(n)})_{t \in [0,T]}$ is due to the KL-divergence term under the integral, since a priori we do not know whether it is integrable. We will now prove by induction how to bound that term. Suppose that $\int_0^T \operatorname{KL}(m_s^{(n-1)}|\pi) ds < \infty$ and, based on \eqref{eq:Picard}, write
		\begin{equation}\label{eq:auxPicard1}
		\begin{split}
		\log \frac{m_t^{(n)}(x)}{\pi(x)} &= e^{-\frac{\sigma^2}{2}t} \log \frac{m_0(x)}{\pi(x)} \\
		&- \int_0^t \frac{\sigma^2}{2} e^{-\frac{\sigma^2}{2}(t-s)} \left( \frac{2}{\sigma^2} \frac{\delta F}{\delta m}(m_s^{(n-1)},x) - \operatorname{KL}(m_s^{(n-1)}|\pi) \right) ds \,.
		\end{split}
		\end{equation}
		We also have
		\begin{equation}\label{eq:auxPicard2}
		\begin{split}
		\log \frac{\pi(x)}{m_t^{(n)}(x)} &= - e^{-\frac{\sigma^2}{2}t} \log \frac{m_0(x)}{\pi(x)} \\
		&- \int_0^t \frac{\sigma^2}{2} e^{-\frac{\sigma^2}{2}(t-s)} \left( - \frac{2}{\sigma^2} \frac{\delta F}{\delta m}(m_s^{(n-1)},x) + \operatorname{KL}(m_s^{(n-1)}|\pi) \right) ds \,.
		\end{split}
		\end{equation}
		Due to \eqref{eq:assumpF} and \eqref{eq:assumpSup}, we can multiply both sides of \eqref{eq:auxPicard1} by $m_t^{(n)}(x)$ and integrate with respect to $x$ in order to obtain
		\begin{equation*}
		\operatorname{KL}(m_t^{(n)}|\pi) \leq \log R + \frac{2}{\sigma^2}C + \int_0^t \frac{\sigma^2}{2} e^{-\frac{\sigma^2}{2}(t-s)}  \operatorname{KL}(m_s^{(n-1)}|\pi) ds \,.
		\end{equation*}
		Similarly, by multiplying both sides of \eqref{eq:auxPicard2} by $\pi(x)$ and integrating with respect to $x$, we obtain
		\begin{equation*}
		\operatorname{KL}(\pi|m_t^{(n)}) \leq \log R + \frac{2}{\sigma^2}C - \int_0^t \frac{\sigma^2}{2} e^{-\frac{\sigma^2}{2}(t-s)}  \operatorname{KL}(m_s^{(n-1)}|\pi) ds \,.
		\end{equation*}
		Consequently, we obtain
		\begin{equation*}
		\operatorname{KL}(m_t^{(n)}|\pi) \leq \operatorname{KL}(m_t^{(n)}|\pi) + \operatorname{KL}(\pi|m_t^{(n)}) \leq 2 \log R + \frac{4}{\sigma^2}C \,,
		\end{equation*}
		which finishes the proof by induction.
	\end{proof}	
\end{lemma}

We will now consider the sequence of flows $\left( (m_t^{(n)})_{t \in [0,T]} \right)_{n=0}^{\infty}$ in $\mathcal{P}(\mathbb{R}^d)^{[0,T]}$ equipped with the distance $\mathcal{TV}_T$, defined for any $(\mu_t)_{t \in [0,T]}$, $(\nu_t)_{t \in [0,T]} \in \mathcal{P}(\mathbb{R}^d)^{[0,T]}$ by
\begin{equation*}
\mathcal{TV}_T \left( (\mu_t)_{t \in [0,T]}, (\nu_t)_{t \in [0,T]} \right) := \int_0^T TV(\mu_t, \nu_t) dt \,.
\end{equation*}
Since $\mathcal{P}(\mathbb{R}^d)$ equipped with the total variation distance $TV$ is complete, we can apply the argument from Lemma A.5 in \cite{SiskaSzpruch2020} with $p=1$ to conclude that $\mathcal{P}(\mathbb{R}^d)^{[0,T]}$ equipped with $\mathcal{TV}_T$ is also complete. We will now consider the Picard iteration mapping $\Psi\left((m_t^{(n-1)})_{t \in [0,T]}\right) := (m_t^{(n)})_{t \in [0,T]}$ defined via \eqref{eq:Picard}, and show that $\Psi$ is contractive in $(\mathcal{P}(\mathbb{R}^d)^{[0,T]},\mathcal{TV}_T)$. Then the Banach fixed point theorem will give us the existence of a solution to \eqref{eq:bd1}.

\begin{lemma}\label{thm:Contractivity}
	The mapping $\Psi\left((m_t^{(n-1)})_{t \in [0,T]}\right) := (m_t^{(n)})_{t \in [0,T]}$ defined via \eqref{eq:Picard} is contractive in $(\mathcal{P}(\mathbb{R}^d)^{[0,T]},\mathcal{TV}_T)$.
	\begin{proof}
		From \eqref{eq:Picard} we have
		\begin{equation*}
		\begin{split}
		&\log m_t^{(n)}(x) - \log m_t^{(n-1)}(x) = - \int_0^t \frac{\sigma^2}{2} e^{-\frac{\sigma^2}{2}(t-s)} \times\\
		&\times \left[  \frac{2}{\sigma^2} \left( \frac{\delta F}{\delta m}(m_s^{(n-1)},x) - \frac{\delta F}{\delta m}(m_s^{(n-2)},x) \right)   - \operatorname{KL}(m_s^{(n-1)}|\pi) + \operatorname{KL}(m_s^{(n-2)}|\pi) \right] ds \,.
		\end{split}
		\end{equation*}
		Multiplying both sides by $m_t^{(n)}(x)$ and integrating with respect to $x$, we obtain
		\begin{equation}\label{eq:contrAux1}
		\begin{split}
		&\operatorname{KL}(m_t^{(n)}|m_t^{(n-1)}) = - \int_0^t \frac{\sigma^2}{2} e^{-\frac{\sigma^2}{2}(t-s)} \Bigg[  \frac{2}{\sigma^2} \int_{\mathbb{R}^d} \left( \frac{\delta F}{\delta m}(m_s^{(n-1)},x) - \frac{\delta F}{\delta m}(m_s^{(n-2)},x) \right)\\
		&\times m_t^{(n)}(dx)    - \operatorname{KL}(m_s^{(n-1)}|\pi) + \operatorname{KL}(m_s^{(n-2)}|\pi) \Bigg] ds \,.
		\end{split}
		\end{equation}
		Moreover, note that 
		\begin{equation*}
		\begin{split}
		&\int_{\mathbb{R}^d} \left( \frac{\delta F}{\delta m}(m_s^{(n-1)},x) - \frac{\delta F}{\delta m}(m_s^{(n-2)},x) \right)m_t^{(n)}(dx)\\
		&= \int_{\mathbb{R}^d} \int_{\mathbb{R}^d} \int_0^1 \frac{\delta^2 F}{\delta m^2} \left( m_s^{(n-2)} + \lambda \left( m_s^{(n-1)} - m_s^{(n-2)} \right),x,y  \right) d \lambda \\ &\times \left(m_s^{(n-1)} - m_s^{(n-2)} \right)(dy) m_t^{(n)}(dx) \,.
		\end{split}
		\end{equation*}
		Similarly, again from \eqref{eq:Picard} we have 
		\begin{equation*}
		\begin{split}
		&\log m_t^{(n-1)}(x) - \log m_t^{(n)}(x) = - \int_0^t \frac{\sigma^2}{2} e^{-\frac{\sigma^2}{2}(t-s)} \times\\
		&\times \left[  \frac{2}{\sigma^2} \left( \frac{\delta F}{\delta m}(m_s^{(n-2)},x) - \frac{\delta F}{\delta m}(m_s^{(n-1)},x) \right)   - \operatorname{KL}(m_s^{(n-2)}|\pi) + \operatorname{KL}(m_s^{(n-1)}|\pi) \right] ds \,.
		\end{split}
		\end{equation*}
		Multiplying both sides by $m_t^{(n-1)}(x)$ and integrating with respect to $x$, we obtain
		\begin{equation}\label{eq:contrAux2}
		\begin{split}
		&\operatorname{KL}(m_t^{(n-1)}|m_t^{(n)}) = - \int_0^t \frac{\sigma^2}{2} e^{-\frac{\sigma^2}{2}(t-s)} \Bigg[  \frac{2}{\sigma^2} \int_{\mathbb{R}^d} \left( \frac{\delta F}{\delta m}(m_s^{(n-2)},x) - \frac{\delta F}{\delta m}(m_s^{(n-1)},x)  \right)\\
		&\times m_t^{(n-1)}(dx)  
		- \operatorname{KL}(m_s^{(n-2)}|\pi) + \operatorname{KL}(m_s^{(n-1)}|\pi) \Bigg] ds \,.
		\end{split}
		\end{equation}
		Similarly as before, we note that
		\begin{equation*}
		\begin{split}
		&\int_{\mathbb{R}^d} \left( \frac{\delta F}{\delta m}(m_s^{(n-2)},x) - \frac{\delta F}{\delta m}(m_s^{(n-1)},x) \right)m_t^{(n-1)}(dx)\\
		&= -\int_{\mathbb{R}^d} \int_{\mathbb{R}^d} \int_0^1 \frac{\delta^2 F}{\delta m^2} \left( m_s^{(n-2)} + \lambda \left( m_s^{(n-1)} - m_s^{(n-2)} \right),x,y  \right) d \lambda\\
		&\times \left(m_s^{(n-1)} - m_s^{(n-2)} \right)(dy) m_t^{(n-1)}(dx) \,.
		\end{split}
		\end{equation*}
		Combining \eqref{eq:contrAux1} and \eqref{eq:contrAux2}, we obtain
		\begin{equation*}
		\begin{split}
		&\operatorname{KL}(m_t^{(n)}|m_t^{(n-1)}) + \operatorname{KL}(m_t^{(n-1)}|m_t^{(n)}) = - \int_0^t e^{-\frac{\sigma^2}{2}(t-s)} \times\\
		&\times \int_{\mathbb{R}^d} \int_{\mathbb{R}^d} \int_0^1 \frac{\delta^2 F}{\delta m^2} \left( m_s^{(n-2)} + \lambda \left( m_s^{(n-1)} - m_s^{(n-2)} \right),x,y  \right) d \lambda \left(m_s^{(n-1)} - m_s^{(n-2)} \right)(dy)\times \\
		&\times \left( m_t^{(n)} - m_t^{(n-1)}\right)(dx)   ds \,.
		\end{split}
		\end{equation*}
		Hence, due to \eqref{eq:assumpF2}, we get
		\begin{equation*}
		\begin{split}
		&\operatorname{KL}(m_t^{(n)}|m_t^{(n-1)}) + \operatorname{KL}(m_t^{(n-1)}|m_t^{(n)}) \\
		&\leq \int_0^t e^{-\frac{\sigma^2}{2}(t-s)} C_2 TV(m_s^{(n-1)},m_s^{(n-2)}) TV(m_t^{(n)},m_t^{(n-1)}) ds
		\end{split}
		\end{equation*}
		By the Pinsker-Csizsar inequality,  $TV^2(m_t^{(n)},m_t^{(n-1)}) \leq \frac{1}{2} \operatorname{KL}(m_t^{(n)}|m_t^{(n-1)})$ and hence
		\begin{equation*}
		4 TV^2(m_t^{(n)},m_t^{(n-1)}) \leq C_2 TV(m_t^{(n)},m_t^{(n-1)}) \int_0^t e^{-\frac{\sigma^2}{2}(t-s)}  TV(m_s^{(n-1)},m_s^{(n-2)}) ds \,,
		\end{equation*}
		which gives
		\begin{equation*}
		\begin{split}
		TV(m_t^{(n)},&m_t^{(n-1)}) \leq \frac{C_2}{4}  \int_0^t e^{-\frac{\sigma^2}{2}(t-s)}  TV(m_s^{(n-1)},m_s^{(n-2)}) ds \\
		&\leq \left( \frac{C_2}{4} \right)^{n-1}  e^{-\frac{\sigma^2}{2}t} \int_0^t \int_0^{t_1} \ldots \int_0^{t_{n-2}}  e^{\frac{\sigma^2}{2}t_{n-1}}  TV(m_{t_{n-1}}^{(1)}, m_{t_{n-1}}^{(0)}) dt_{n-1} \ldots dt_2 dt_1 \\
		&\leq \left( \frac{C_2}{4} \right)^{n-1}  e^{-\frac{\sigma^2}{2}t} \frac{t^{n-2}}{(n-2)!}\int_0^t   e^{\frac{\sigma^2}{2}t_{n-1}}  TV(m_{t_{n-1}}^{(1)}, m_{t_{n-1}}^{(0)}) dt_{n-1} \\
		&\leq \left( \frac{C_2}{4} \right)^{n-1}   \frac{t^{n-2}}{(n-2)!}\int_0^t    TV(m_{t_{n-1}}^{(1)}, m_{t_{n-1}}^{(0)}) dt_{n-1} \,,
		\end{split}
		\end{equation*}
		where in the third inequality we bounded $\int_0^{t_{n-2}}  dt_{n-1} \leq \int_0^{t}  dt_{n-1}$ and in the fourth inequality we bounded $e^{\frac{\sigma^2}{2}t_{n-1}}  \leq e^{\frac{\sigma^2}{2}t}$. Hence we obtain
		\begin{equation*}
		\int_0^T	TV(m_t^{(n)},m_t^{(n-1)})dt \leq \left( \frac{C_2}{4} \right)^{n-1}   \frac{T^{n-1}}{(n-2)!}\int_0^T    TV(m_{t_{n-1}}^{(1)}, m_{t_{n-1}}^{(0)}) dt_{n-1} \,.
		\end{equation*}
		For sufficiently large $n$, the constant on the right hand side becomes less than $1$ and the proof is complete.
	\end{proof}	
\end{lemma} 

We can now finalize the proof of Theorem \ref{thm:existence}.

\begin{proofExistence} 
	\emph{Step 1: Existence of the gradient flow and bound \eqref{eq:KLmtBound} on $[0,T]$.}
		By Lemma \ref{thm:Contractivity}, for any $T > 0$ we obtain 
	the existence of a flow $(m_t)_{t \in [0,T]}$ satisfying \eqref{eq:bd1}. Moreover, for Lebesgue-almost all $t \in [0,T]$ we have
	\begin{equation*}
	TV(m_t^{(n)},m_t) \to 0 \qquad \text{ as } n \to \infty \,,
	\end{equation*}
	which implies
	\begin{equation*}
	m_t^{(n)} \to m_t \qquad \text{ weakly, as } n \to \infty \,.
	\end{equation*}
	Hence, using the lower semi-continuity of the KL-divergence (see e.g.\ Theorem 2.34 in \cite{AmbrosioFuscoPallara2000}) we obtain
	\begin{equation}\label{eq:KLmtBound2}
	\operatorname{KL}(m_t|\pi) \leq \liminf_{n \to \infty}  \operatorname{KL}(m_t^{(n)}|\pi) \leq 2 \log R + \frac{4C}{\sigma^2} \,, 
	\end{equation}
	where the second inequality follows from Lemma \ref{lem:existenceGF}. In order to ensure that the solution $(m_t)_{t \in [0,T]}$ can be extended to all $t \geq 0$, we first need to prove the bound on the ratio $m_t/\pi$ in \eqref{eq:ratioBoundAllt}.

	\emph{Step 2: Ratio condition \eqref{eq:ratioBoundAllt}.}
			Following the discussion from the beginning of Section \ref{section:existence}, we see that for any $t \in [0,T]$ we have 
			\begin{equation*}
			\begin{split}
			\log \frac{m_t(x)}{\pi(x)} &= e^{-\frac{\sigma^2}{2}t} \log \frac{m_0(x)}{\pi(x)} \\
			&- \int_0^t \frac{\sigma^2}{2} e^{-\frac{\sigma^2}{2}(t-s)} \left( \frac{2}{\sigma^2} \frac{\delta F}{\delta m}(m_s,x) - \operatorname{KL}(m_s|\pi) \right) ds \,.
			\end{split}
			\end{equation*}
			Using \eqref{eq:assumpF}, \eqref{eq:assumpSup} and \eqref{eq:KLmtBound2} we obtain
			\begin{equation*}
			\log \frac{m_t(x)}{\pi(x)} \leq \log R + C + \frac{\sigma^2}{2} \left( 2 \log R + \frac{4C}{\sigma^2} \right) \,. 
			\end{equation*}
			Hence we can choose $R_1 := 1 + \exp \left( \log R + C + \frac{\sigma^2}{2} \left( 2 \log R + \frac{4C}{\sigma^2} \right) \right)$. Note that we choose $R_1 > 1$ purely for convenience, to ensure that $\log R_1 > 0$ in our subsequent calculations.  
			Obtaining a lower bound on $\frac{m_t(x)}{\pi(x)}$ follows similarly, by using \eqref{eq:assumpInf} instead of \eqref{eq:assumpSup}.
	
\emph{Step 3: Existence of the gradient flow on $[0,\infty)$.}	
 In order to complete our proof, note that the unique solution $(m_t)_{t \in [0,T]}$ to \eqref{eq:bd1} can also be expressed as
	\begin{equation*}
	m_t(x) = m_0(x) \exp \left( - \int_0^t \left( \frac{\delta F}{\delta m} (m_s,x) + \frac{\sigma^2}{2}\log \left( \frac{m_s(x)}{\pi(x)}\right) - \frac{\sigma^2}{2}\operatorname{KL}(m_s|\pi) \right) ds \right) \,. 
	\end{equation*}
	From \eqref{eq:assumpF}, \eqref{eq:KLmtBound2} and \eqref{eq:ratioBoundAllt}, we obtain for any $t \in [0,T]$
	\begin{equation*}
	\begin{split}
	&\left| \frac{\delta F}{\delta m}(m_t,x) + \frac{\sigma^2}{2}\log \left( \frac{m_t(x)}{\pi(x)} \right) - \frac{\sigma^2}{2}\operatorname{KL}(m_t|\pi) \right| \\
	&\leq 3C + \frac{\sigma^2}{2}\left( \max\{ |\log r_1|,\log R_1 \} + 2 \log R \right)  =: C_V \,.
	\end{split}
	\end{equation*}
	This gives $\| m_t \|_{TV} \leq \| m_0 \|_{TV} e^{C_V t}$, and shows that $m_t$ does not explode in any finite time, hence we obtain a global solution $(m_t)_{t \in [0,\infty)}$. In particular, the bounds in \eqref{eq:KLmtBound2}, \eqref{eq:ratioBoundAllt} and \eqref{eq:ratioBoundAllt2} hold for all $t \geq 0$.
	\qed
\end{proofExistence}

\begin{proofDifferentiability}
	We have the differentiability of $F(m_t)$ as a consequence of Assumption \ref{assump:F}. 
	In order to show the differentiability of $\operatorname{KL}(m_t|\pi) = \int_{\mathbb{R}^d} \log \left( \frac{m_t(x)}{\pi(x)}\right) m_t(x) dx = \int_{\mathbb{R}^d} \log \left( \frac{m_t(x)}{\pi(x)}\right) \frac{m_t(x)}{\pi(x)} \pi(x)dx$, we need to prove that
	$\left| \partial_t \left( \log \left( \frac{m_t(x)}{\pi(x)}\right) \frac{m_t(x)}{\pi(x)}\right) \right| \leq g(x)$
	for some function $g$ integrable with respect to $\pi$, which is sufficient 
	by a standard result in measure theory (see e.g.\ Theorem 11.5 in \cite{Schilling2005}).
	Indeed, by \eqref{eq:assumpF}, \eqref{eq:ratioBoundAllt} and \eqref{eq:KLmtBound}, we get
	\begin{equation*}
	\begin{split}
	&\left| \partial_t \left( \log \left( \frac{m_t(x)}{\pi(x)}\right) \frac{m_t(x)}{\pi(x)} \right) \right| = \left| \frac{\pi(x)}{m_t(x)}\frac{\partial_t m_t(x)}{\pi(x)} \frac{m_t(x)}{\pi(x)} + \log \left( \frac{m_t(x)}{\pi(x)}\right) \frac{\partial_t m_t(x)}{\pi(x)} \right| \\
	&= \left| \left( 1+\log\left(\frac{m_t(x)}{\pi(x)}\right)\right) \frac{\partial_t m_t(x)}{\pi(x)} \right| \\
	&= \left| \left( 1 + \log \left( \frac{m_t(x)}{\pi(x)} \right) \right) \left( \frac{\delta F}{\delta m}(m_t,x) + \frac{\sigma^2}{2}\log \left( \frac{m_t(x)}{\pi(x)} \right) - \frac{\sigma^2}{2}\operatorname{KL}(m_t|\pi)\right) \frac{m_t(x)}{\pi(x)} \right|\\
	&\leq \left( 1 + \max\{ |\log r_1|, \log R_1 \} \right) \left( 3C + \frac{\sigma^2}{2}\left( \max\{ |\log r_1|, \log R_1 \} + 2 \log R \right) \right) R_1 \,. 
	\end{split}
	\end{equation*}
	We can now write
	\begin{equation*}
	\begin{split}
	\partial_t V^{\sigma}(m_t) &= \int_{\mathbb{R}^d} \frac{\delta F}{\delta m}(m_t,x) \partial_t m_t(x) dx + \frac{\sigma^2}{2}\int_{\mathbb{R}^d} \partial_t \left( \log \left( \frac{m_t(x)}{\pi(x)} \right) m_t(x) \right) dx \\
	&= \int_{\mathbb{R}^d} \left[ \frac{\delta F}{\delta m}(m_t,x) +  \frac{\sigma^2}{2}\left( 1 + \log \left( \frac{m_t(x)}{\pi(x)}\right) \right) \right] \partial_t m_t(x) dx \\
	&= \int_{\mathbb{R}^d} \left( \frac{\delta F}{\delta m} (m_t,x) + \frac{\sigma^2}{2} \log \left( \frac{m_t(x)}{\pi(x)}\right) - \frac{\sigma^2}{2} \operatorname{KL}(m_t|\pi) \right) \partial_t m_t(x) dx 
	\,,
	\end{split}
	\end{equation*}
	where the last equality follows due to $\int_{\mathbb{R}^d} \partial_t m_t(x) dx = 0$. Combining this with \eqref{eq:birthdeath} proves  \eqref{eq:diff}. 
	\qed
\end{proofDifferentiability}

\begin{proofPLI}
	By Lemma \ref{lem:quadraticGrowth}, the quadratic growth condition \eqref{eq:quadraticGrowthKL} required in Theorem \ref{thm:mainLojasiewicz} is satisfied for $m=m_t$ for all $t >0$, with $\lambda = \sigma^2/2$. Moreover, due to \eqref{eq:ratioBoundAllt} and \eqref{eq:ratioBoundAllt2}, the ratio condition \eqref{eq:ratio} required in Theorem \ref{thm:mainLojasiewicz} is satisfied for $m=m_t$ for all $t >0$. Indeed, recall that by the discussion below Assumption \ref{assump:m0}, the ratio condition for $m_0/\pi$ with constants $r$ and $R$ is equivalent to the ratio condition for $m_0/m^{\sigma,*}$ with corresponding constants $\bar{r}$ and $\bar{R}$. Similarly, since due to \eqref{eq:ratioBoundAllt} and \eqref{eq:ratioBoundAllt2} we have a bound on $m_t/\pi$ for all $t >0$ with constants $r_1$ and $R_1$, we can apply the argument below Assumption \ref{assump:m0} to obtain a bound on $m_t/m^{\sigma,*}$ for all $t >0$, with appropriately modified constants $\bar{r}_1$ and $\bar{R}_1$. Furthermore, note that by the proof of Lemma \ref{lem:quadraticGrowth}, in the case of $G=V^{\sigma}$, the convexity condition needed in the proof of Lemma \ref{lemma:flatCauchySchwarz} (and thus in Theorem \ref{thm:mainLojasiewicz}) can be applied with $a(m,x)$ instead of $\frac{\delta G}{\delta m}$, i.e., for any $m$, $m' \in \mathcal{P}(\mathbb{R}^d)$ we have
	\begin{equation*}
	V^{\sigma}(m) - V^{\sigma}(m') \leq \int_{\mathbb{R}^d} a(m,x)(m-m')(dx) \,.
	\end{equation*}
	 As a consequence, the argument from the proof of Theorem \ref{thm:mainLojasiewicz} applies to our setting and the flat Polyak-\L ojasiewicz condition \eqref{eq:flatLojasiewicz2} is satisfied for all $t \geq 0$.
	\qed
\end{proofPLI}

\section{Appendix: Relations between different f-divergences}\label{section:appendix1}

Suppose we have absolutely continuous probability measures $m$, $m' \in \mathcal{P}(\mathbb{R}^d)$ and a convex function $f: [0,\infty) \to \mathbb{R}$. Then the $f$-divergence of $m$ with respect to $m'$ is defined by
\begin{equation*}
I_f(m|m') := \int_{\mathbb{R}^d} f\left( \frac{m(x)}{m'(x)} \right) m'(x) dx \,. 
\end{equation*}
For instance, choosing $f(t) = t \log t$ leads to the KL-divergence and $f(t) = (t-1)^2$ leads to the $\chi^2$-divergence. We have the following result adapted from Theorem 6 in \cite{Dragomir}.

\begin{lemma}\label{lem:DragomirAux}
	Let $f: [0,\infty) \to \mathbb{R}$ be convex and such that $f(1)=0$. Let us consider an interval $(r,R) \subset (0,\infty)$ such that\\
	(i) $f$ is twice differentiable on $(r,R)$\\
	(ii) there exist real constants $a$, $A$ such that
	$$a \leq t f''(t) \leq A \quad \text{for all } t \in (r,R) \,.$$
	Then for any absolutely continuous probability measures $\mu$ and $\nu$, we have the inequality 
	$$a \operatorname{KL}(\mu|\nu) \leq I_f(\mu|\nu) \leq A  \operatorname{KL}(\mu|\nu) \,.$$
	\begin{proof}
		Let us define a mapping $F_a: (0,\infty) \to \mathbb{R}$ given by $F_a(t) := f(t)-a t\log t$. Then $F_a$ is such that $F_a(1)=0$, and is twice differentiable and convex on $(r,R)$ since $F_a''(t) \geq 0$ on $(r,R)$. 
		Note that the $f$-divergence associated to a convex $F_a$ with $F_a(1)=0$ is always non-negative due to Jensen's inequality, and hence 
		we have
		$$0 \leq I_{F_a}(\mu|\nu) = I_f(\mu|\nu) - a \operatorname{KL} (\mu|\nu) \,.$$
		We now define a mapping $F_A: (0,\infty) \to \mathbb{R}$ by setting $F_A(t) := A t\log t -f(t)$. Then $F_A$ is such that $F_A(1)=0$, and is twice differentiable and convex on $(r,R)$ since $F_A''(t) \geq 0$ on $(r,R)$. We again use the fact that the corresponding $f$-divergence is non-negative, and we obtain
		$$0 \leq I_{F_A}(\mu|\nu) = A \operatorname{KL} (\mu|\nu)-I_f(\mu|\nu)\,,$$
		which finishes the proof.
	\end{proof}
\end{lemma}

\begin{proofDragomir}
	We consider the mapping $f_1: (0,\infty) \to \mathbb{R}$ given by $f_1(t)= -\log(t)$. Note that the $f$-divergence corresponding to this $f_1$ is the KL-divergence with swapped arguments, i.e., for any absolutely continuous probability measures $\mu$ and $\nu$, we have
	\begin{equation*}
	I_{f_1}(\mu|\nu) = -\int_{\mathbb{R}^d} \log \frac{\mu(x)}{\nu(x)} \nu(x) dx = \int_{\mathbb{R}^d} \log \frac{\nu(x)}{\mu(x)} \nu(x) dx = \operatorname{KL}(\nu|\mu) \,.
	\end{equation*}
	We remark that $f_1(1)=0$ and that $f_1$ is twice differentiable on any interval $(r,R) \subset (0,\infty)$. We also have $$\frac{1}{R} \leq t f_1''(t) \leq \frac{1}{r} \quad \text{for all } t \in (r,R)$$ since $f''(t)=1/t^2$. Applying Lemma \ref{lem:DragomirAux} with $a = 1/R$ and $A = 1/r$, we have $$\frac{1}{R} \operatorname{KL}(\mu|\nu) \leq \operatorname{KL} (\nu|\mu) \leq \frac{1}{r}  \operatorname{KL}(\mu|\nu) \,.$$
	This shows the first inequality in \eqref{eq:Dragomir}. We now consider the mapping $f_2: (0,\infty) \to \mathbb{R}$ defined by $f_2(t):= (t-1)^2$, i.e., $I_{f_2}$ is the $\chi^2$-divergence. Again, $f_2(1)=0$ and $f_2$ is twice differentiable on any interval $(r,R) \subset (0,\infty)$. Moreover, we have 
	$$2r \leq t f_2''(t) \leq 2R \quad \text{for all } t \in(r,R)$$ since $f_2''=2$. Applying Lemma \ref{lem:DragomirAux} with $a=2r$ and $A=2R$, we have 
	$$2r \operatorname{KL}(\mu|\nu) \leq \chi^2 (\mu|\nu) \leq 2 R  \operatorname{KL}(\mu|\nu) \,. $$
	This shows the second inequality in \eqref{eq:Dragomir} and concludes the proof.
	\qed
\end{proofDragomir}

\section{Appendix: Flat Derivative}\label{section:appendix2}

\begin{definition}\label{def fderivative} Fix $q\geq 0$ and let $\mathcal P_q(\mathbb{R}^d)$ be the space of probability measures on $\mathbb{R}^d$ with finite $q$-th moments. A functional $F:\mathcal P_q(\mathbb{R}^d) \to \mathbb R$, is said to admit a first order linear derivative (or a flat derivative), if there exists a functional $\frac{\delta F}{\delta m}: \mathcal P_q(\mathbb{R}^d) \times\mathbb R^d\rightarrow \mathbb R$, such that
\begin{enumerate}
\item For all $a\in \mathbb R^d$, $\mathcal  P_q(\mathbb{R}^d) \ni m \mapsto \frac{\delta F}{\delta m}(m,a)$ is continuous.
\item For any $\nu \in \mathcal P_q(\mathbb{R}^d)$, there exists $C>0$ such that for all $a\in \mathbb R^d$ we have 
\[
\left|\frac{\delta F}{\delta m}({\nu},a)\right|\leq C(1+|a|^q)\,.
\]
\item For all $m$, $m'\in \mathcal P_q (\mathbb{R}^d)$,
\begin{equation}\label{def FlatDerivative}
F(m')-F(m)=\int_{0}^{1}\int_{\mathbb{R}^d}\frac{\delta F}{\delta m}(m + \lambda (m'-m),a)\left(m'- m\right)(da)\,d\lambda.
\end{equation}
\end{enumerate}
The functional $\frac{\delta F}{\delta m}$ is then called the linear (functional) derivative of $F$ on $\mathcal P_q(\mathbb{R}^d)$.
\end{definition}

Note that Definition \ref{def fderivative} easily generalizes to higher order linear derivatives. More precisely, for a fixed $a \in \mathbb{R}^d$ the functional $\frac{\delta F}{\delta m}(\cdot, a): \mathcal P_q(\mathbb{R}^d) \rightarrow \mathbb R$ can admit a first order linear derivative $\frac{\delta}{\delta m}\left(\frac{\delta F}{\delta m}(\cdot, a)\right): \mathcal P_q(\mathbb{R}^d) \times \mathbb{R}^d \rightarrow \mathbb R$ whenever the conditions from Definition \ref{def fderivative} are satisfied. If that derivative exists for any $a \in \mathbb{R}^d$, we say that $F$ admits a second order linear derivative $\frac{\delta^2 F}{\delta m^2}: \mathcal P_q(\mathbb{R}^d) \times \mathbb{R}^d \times \mathbb{R}^d \rightarrow \mathbb R$, which then satisfies, for all $a \in \mathbb{R}^d$, and for all $m$, $m'\in \mathcal P_q (\mathbb{R}^d)$
\begin{equation*}
\frac{\delta F}{\delta m}(m',a)-\frac{\delta F}{\delta m}(m,a)=\int_{0}^{1}\int_{\mathbb{R}^d}\frac{\delta^2 F}{\delta m^2}(m + \lambda (m'-m),a',a)\left(m'- m\right)(da')\,d\lambda.
\end{equation*}

\section*{Acknowledgements}
LS acknowledges the support of the UKRI Prosperity Partnership Scheme (FAIR) under EPSRC Grant EP/V056883/1.
 
\section*{Declarations}

\subsection*{Competing Interests} The authors have no relevant financial or non-financial interests to disclose.

\bibliographystyle{abbrv}
\bibliography{birthdeath}

\def\cprime{$'$}
\begin{thebibliography}{10}

\bibitem{AmbrosioFuscoPallara2000}
L.~Ambrosio, N.~Fusco, and D.~Pallara.
\newblock {\em Functions of bounded variation and free discontinuity problems}.
\newblock Oxford Mathematical Monographs. The Clarendon Press, Oxford
  University Press, New York, 2000.

\bibitem{AmbrosioGigliSavare2008}
L.~Ambrosio, N.~Gigli, and G.~Savar\'{e}.
\newblock {\em Gradient flows in metric spaces and in the space of probability
  measures}.
\newblock Lectures in Mathematics ETH Z\"{u}rich. Birkh\"{a}user Verlag, Basel,
  second edition, 2008.

\bibitem{Aubin-Frankowski2022}
P.-C. {Aubin-Frankowski}, A.~{Korba}, and F.~{L{\'e}ger}.
\newblock {Mirror Descent with Relative Smoothness in Measure Spaces, with
  application to Sinkhorn and EM}.
\newblock {\em arXiv e-prints}, page arXiv:2206.08873, June 2022.

\bibitem{BenamouBrenier2000}
J.-D. Benamou and Y.~Brenier.
\newblock A computational fluid mechanics solution to the {M}onge-{K}antorovich
  mass transfer problem.
\newblock {\em Numer. Math.}, 84(3):375--393, 2000.

\bibitem{BlanchetBolte2018}
A.~Blanchet and J.~Bolte.
\newblock A family of functional inequalities: {\L}ojasiewicz inequalities and
  displacement convex functions.
\newblock {\em J. Funct. Anal.}, 275(7):1650--1673, 2018.

\bibitem{Bolte2010}
J.~Bolte, A.~Daniilidis, O.~Ley, and L.~Mazet.
\newblock Characterizations of {\l}ojasiewicz inequalities: subgradient flows,
  talweg, convexity.
\newblock {\em Trans. Amer. Math. Soc.}, 362(6):3319--3363, 2010.

\bibitem{CarmonaDelarue2018}
R.~Carmona and F.~Delarue.
\newblock {\em Probabilistic theory of mean field games with applications.
  {I}}, volume~83 of {\em Probability Theory and Stochastic Modelling}.
\newblock Springer, Cham, 2018.
\newblock Mean field FBSDEs, control, and games.

\bibitem{Chizat2021}
L.~{Chizat}.
\newblock {Convergence Rates of Gradient Methods for Convex Optimization in the
  Space of Measures}.
\newblock {\em arXiv e-prints}, page arXiv:2105.08368, May 2021.

\bibitem{Chizat2022}
L.~{Chizat}.
\newblock {Mean-Field Langevin Dynamics: Exponential Convergence and
  Annealing}.
\newblock {\em arXiv e-prints}, page arXiv:2202.01009, Feb. 2022.

\bibitem{ChizatBach}
L.~Chizat and F.~Bach.
\newblock On the global convergence of gradient descent for over-parameterized
  models using optimal transport.
\newblock In S.~Bengio, H.~Wallach, H.~Larochelle, K.~Grauman, N.~Cesa-Bianchi,
  and R.~Garnett, editors, {\em Advances in Neural Information Processing
  Systems}, volume~31. Curran Associates, Inc., 2018.

\bibitem{Dragomir}
S.~{Dragomir}.
\newblock {Upper and lower bounds for Csiszar's f-divergence in terms of the
  Kullback-Leibler distance and applications}.
\newblock 1999.

\bibitem{Gallouet2017}
T.~O. Gallou\"{e}t and L.~Monsaingeon.
\newblock A {JKO} splitting scheme for {K}antorovich-{F}isher-{R}ao gradient
  flows.
\newblock {\em SIAM J. Math. Anal.}, 49(2):1100--1130, 2017.

\bibitem{HuRenSiskaSzpruch2021}
K.~Hu, Z.~Ren, D.~\v{S}i\v{s}ka, and L.~Szpruch.
\newblock Mean-field {L}angevin dynamics and energy landscape of neural
  networks.
\newblock {\em Ann. Inst. Henri Poincar\'{e} Probab. Stat.}, 57(4):2043--2065,
  2021.

\bibitem{Karimi2016}
H.~Karimi, J.~Nutini, and M.~Schmidt.
\newblock Linear convergence of gradient and proximal-gradient methods under
  the {P}olyak-{\l}ojasiewicz condition.
\newblock In P.~Frasconi, N.~Landwehr, G.~Manco, and J.~Vreeken, editors, {\em
  Machine Learning and Knowledge Discovery in Databases}, pages 795--811, Cham,
  2016. Springer International Publishing.

\bibitem{Liero2018}
M.~Liero, A.~Mielke, and G.~Savar\'{e}.
\newblock Optimal entropy-transport problems and a new
  {H}ellinger-{K}antorovich distance between positive measures.
\newblock {\em Invent. Math.}, 211(3):969--1117, 2018.

\bibitem{LuLuNolen2019}
Y.~{Lu}, J.~{Lu}, and J.~{Nolen}.
\newblock {Accelerating Langevin Sampling with Birth-death}.
\newblock {\em arXiv e-prints}, page arXiv:1905.09863, May 2019.

\bibitem{Montanari}
S.~Mei, A.~Montanari, and P.-M. Nguyen.
\newblock A mean field view of the landscape of two-layer neural networks.
\newblock {\em Proc. Natl. Acad. Sci. USA}, 115(33):E7665--E7671, 2018.

\bibitem{Nitanda2022}
A.~{Nitanda}, D.~{Wu}, and T.~{Suzuki}.
\newblock {Convex Analysis of the Mean Field Langevin Dynamics}.
\newblock {\em arXiv e-prints}, page arXiv:2201.10469, Jan. 2022.

\bibitem{Radhakrishnan2020}
A.~{Radhakrishnan}, M.~{Belkin}, and C.~{Uhler}.
\newblock {Linear Convergence of Generalized Mirror Descent with Time-Dependent
  Mirrors}.
\newblock {\em arXiv e-prints}, page arXiv:2009.08574, Sept. 2020.

\bibitem{RenWang2022}
Z.~{Ren} and S.~{Wang}.
\newblock {Entropic fictitious play for mean field optimization problem}.
\newblock {\em arXiv e-prints}, page arXiv:2202.05841, Feb. 2022.

\bibitem{vandenEijnden2019}
G.~Rotskoff, S.~Jelassi, J.~Bruna, and E.~Vanden-Eijnden.
\newblock Neuron birth-death dynamics accelerates gradient descent and
  converges asymptotically.
\newblock In K.~Chaudhuri and R.~Salakhutdinov, editors, {\em Proceedings of
  the 36th International Conference on Machine Learning}, volume~97 of {\em
  Proceedings of Machine Learning Research}, pages 5508--5517. PMLR, 09--15 Jun
  2019.

\bibitem{santambrogio2017euclidean}
F.~Santambrogio.
\newblock $\{$Euclidean, metric, and Wasserstein$\}$ gradient flows: an
  overview.
\newblock {\em Bulletin of Mathematical Sciences}, 7(1):87--154, 2017.

\bibitem{Schilling2005}
R.~L. Schilling.
\newblock {\em Measures, integrals and martingales}.
\newblock Cambridge University Press, New York, 2005.

\bibitem{Villani}
C.~Villani.
\newblock {\em Optimal transport}, volume 338 of {\em Grundlehren der
  mathematischen Wissenschaften [Fundamental Principles of Mathematical
  Sciences]}.
\newblock Springer-Verlag, Berlin, 2009.
\newblock Old and new.

\bibitem{SiskaSzpruch2020}
D.~{{\v{S}}i{\v{s}}ka} and {\L}.~{Szpruch}.
\newblock {Gradient Flows for Regularized Stochastic Control Problems}.
\newblock {\em arXiv e-prints}, page arXiv:2006.05956, June 2020.

\end{thebibliography}

\end{document}